\documentclass[12pt,psamsfonts,leqno,oneside,letterpaper]{amsart}
\usepackage[dvips,text={6.5truein,9truein},left=1truein,top=1truein]{geometry}
\usepackage{amssymb,amsmath,amscd,enumerate,
}
\usepackage[pdftex]{graphicx}
\usepackage{url}

\usepackage[colorlinks,linkcolor=blue,citecolor=blue,pdfstartview=FitH]{hyperref}
\input xy
\xyoption{all}
\SelectTips{cm}{12}

\usepackage{color}
\usepackage{mathrsfs}
\newcommand\redden[1]{{\color{red}#1}}

\newcommand\invisiblecomment[1]{\empty}

\newcommand\missingref[1]{\empty}

\newcommand\tinymissingref[1]{\empty}
\newcommand\abstractcomment[1]{\empty}

\theoremstyle{definition}
\newtheorem{para}{}[section]

\newtheorem{remark}[para]{Remark}
\newtheorem*{wilfrei}{Wilson's Freiheitssatz}
\newtheorem{reformulation}[para]{Reformulation}
\newtheorem{remarks}[para]{Remarks}
\newtheorem{notation}[para]{Notation}

\newtheorem{convention}[para]{Convention}
\newtheorem{definition}[para]{Definition}
\newtheorem{definitions}[para]{Definitions}
\newtheorem{definitionnotation}[para]{Definition and Notation}
\newtheorem{notationdefinitionremark}[para]{Notation, Definition and Remark}

\newtheorem{remarksnotation}[para]{Remarks and Notation}

\newtheorem{remarknotation}[para]{Remark and Notation}
\newtheorem{notationremark}[para]{Notation and Remark}
\newtheorem{notationreviewremarks}[para]{Notation, Review and Remarks}
\newtheorem{definitionremark}[para]{Definition and Remark}
\newtheorem{definitionsremarks}[para]{Definitions and Remarks}
\newtheorem{notationremarks}[para]{Notation and Remarks}
\newtheorem{definitionsnotation}[para]{Definitions and Notation}

\newtheorem{reviewdefinition}[para]{Review and Definition}
\newtheorem{definitionremarks}[para]{Definition and Remarks}
\newtheorem{definitionnotationremarks}[para]{Definition, Notation and Remarks}
\newtheorem{definitionsnotationremarks}[para]{Definitions, Notation and Remarks}

\newcommand\Alternatives{\begin{enumerate}[(i)]}
\newcommand\EndAlternatives{\end{enumerate}}
\newcommand\Conditions{\begin{enumerate}[(1)]}
\newcommand\EndConditions{\end{enumerate}}

\theoremstyle{plain}
\newtheorem{theorem}[para]{Theorem}

\newtheorem{lemma}[para]{Lemma}
\newtheorem{remarkdefinition}[para]{Remark and Definition}
\newtheorem{remarkdefinitions}[para]{Remark and Definitions}
\newtheorem{proposition}[para]{Proposition}

\newtheorem{corollary}[para]{Corollary}
\newtheorem{conjecture}[para]{Conjecture}

\newtheorem*{theoremA}{Theorem A}
\newtheorem*{theoremB}{Theorem B}

\newtheorem{claim}[equation]{}
\numberwithin{equation}{para}
\numberwithin{figure}{section}
\numberwithin{specialremark}{para}
\numberwithin{specialnumber}{para}

\newcommand\Number{\begin{para}}
\newcommand\EndNumber{\end{para}}
\newcommand\Definition{\begin{definition}}
\newcommand\EndDefinition{\end{definition}}
\newcommand\Definitions{\begin{definitions}}
\newcommand\DefinitionsNotation{\begin{definitionsnotation}}
\newcommand\NotationDefinitionRemark{\begin{notationdefinitionremark}}
\newcommand\EndNotationDefinitionRemark{\end{notationdefinitionremark}}

\newcommand\DefinitionNotation{\begin{definitionnotation}}
\newcommand\RemarksNotation{\begin{remarksnotation}}
\newcommand\Reformulation{\begin{reformulation}}
\newcommand\EndRemarksNotation{\end{remarksnotation}}
\newcommand\EndReformulation{\end{reformulation}}
\newcommand\RemarkNotation{\begin{remarknotation}}
\newcommand\EndRemarkNotation{\end{remarknotation}}
\newcommand\NotationRemark{\begin{notationremark}}
\newcommand\EndDefinitionNotationRemarks{\end{definitionnotationremarks}}
\newcommand\NotationReviewRemarks{\begin{notationreviewremarks}}
\newcommand\DefinitionRemark{\begin{definitionremark}}
\newcommand\DefinitionRemarks{\begin{definitionremarks}}    
\newcommand\DefinitionsRemarks{\begin{definitionsremarks}}
\newcommand\DefinitionNotationRemarks{\begin{definitionnotationremarks}}
\newcommand\DefinitionsNotationRemarks{\begin{definitionsnotationremarks}}
\newcommand\EndDefinitionsNotationRemarks{\end{definitionsnotationremarks}}

\newcommand\NotationRemarks{\begin{notationremarks}}
\newcommand\EndNotationRemark{\end{notationremark}}
\newcommand\EndNotationReviewRemarks{\end{notationreviewremarks}}
\newcommand\EndDefinitionRemark{\end{definitionremark}}
\newcommand\EndDefinitionRemarks{\end{definitionremarks}}
\newcommand\EndDefinitionsRemarks{\end{definitionsremarks}}
\newcommand\EndNotationRemarks{\end{notationremarks}}
\newcommand\EndRemarkDefinition{\end{remarkdefinition}}
\newcommand\EndRemarkDefinitions{\end{remarkdefinitions}}
\newcommand\RemarkDefinition{\begin{remarkdefinition}}
\newcommand\RemarkDefinitions{\begin{remarkdefinitions}}
\newcommand\EndDefinitionsNotation{\end{definitionsnotation}}
\newcommand\EndDefinitionNotation{\end{definitionnotation}}
\newcommand\ReviewDefinition{\begin{reviewdefinition}}
\newcommand\EndReviewDefinition{\end{reviewdefinition}}
\newcommand\EndDefinitions{\end{definitions}}
\newcommand\Theorem{\begin{theorem}}
\newcommand\EndTheorem{\end{theorem}}
\newcommand\Conjecture{\begin{conjecture}}
\newcommand\EndConjecture{\end{conjecture}}
\newcommand\Remark{\begin{remark}}
\newcommand\EndRemark{\end{remark}}
\newcommand\Remarks{\begin{remarks}}
\newcommand\EndRemarks{\end{remarks}}
\newcommand\Convention{\begin{convention}}
\newcommand\EndConvention{\end{convention}}
\newcommand\Notation{\begin{notation}}
\newcommand\EndNotation{\end{notation}}
\newcommand\Lemma{\begin{lemma}}
\newcommand\EndLemma{\end{lemma}}
\newcommand\Proposition{\begin{proposition}}
\newcommand\EndProposition{\end{proposition}}
\newcommand\Corollary{\begin{corollary}}
\newcommand\EndCorollary{\end{corollary}}
\newcommand\Claim{\begin{claim}}
\newcommand\EndClaim{\end{claim}}
\newcommand\Proof{\begin{proof}}
\newcommand\EndProof{\end{proof}}
\newcommand\Equation{\begin{equation}}
\newcommand\EndEquation{\end{equation}}

\newcommand\Bullets{\begin{itemize}}
\newcommand\EndBullets{\end{itemize}}

\renewcommand\epsilon{\varepsilon}

\newcommand\defish{\mathop{\rm def}}

\newcommand\discup{\mathbin{\rotatebox[origin=c]{90}{$\vDash$}}}

\newcommand\chibar{\overline\chi}

\newcommand\barX{\overline{X}}

\newcommand\tM{\widetilde M}

\newcommand\barx{\overline{x}}

\newcommand\barrho{\overline{\rho}}
\newcommand\barsigma{\overline{\sigma}}

\newcommand\barR{{\overline{R}}}
\newcommand\bart{{\overline{t}}}

\newcommand\ZZ{{\bf Z}}

\newcommand\CC{{\mathbb C}}
\newcommand\QQ{{\mathbb Q }}
\newcommand\HH{{\bf H}}

\newcommand\dist{\mathop{\rm dist}}

\newcommand\length{\mathop{{\rm length}}}

\newcommand\pizzle{{\rm PSL}_2}
\newcommand\zzle{{\rm SL}_2}
\newcommand\ggle{{\rm GL}_2}

\newcommand\iof{{\rm iof}}
\newcommand\miof{{\rm miof}}

\newcommand\g{g}

\DeclareFontFamily{U}{rcjhbltx}{}
\DeclareFontShape{U}{rcjhbltx}{m}{n}{<->rcjhbltx}{}
\DeclareSymbolFont{hebrewletters}{U}{rcjhbltx}{m}{n}

\let\aleph\relax\let\beth\relax
\let\gimel\relax\let\daleth\relax

\DeclareMathSymbol{\aleph}{\mathord}{hebrewletters}{39}
\DeclareMathSymbol{\beth}{\mathord}{hebrewletters}{98}
\DeclareMathSymbol{\gimel}{\mathord}{hebrewletters}{103}
\DeclareMathSymbol{\daleth}{\mathord}{hebrewletters}{100}

\DeclareMathSymbol{\lamed}{\mathord}{hebrewletters}{108}
\DeclareMathSymbol{\mem}{\mathord}{hebrewletters}{109}
\DeclareMathSymbol{\ayin}{\mathord}{hebrewletters}{96}
\DeclareMathSymbol{\tsadi}{\mathord}{hebrewletters}{118}
\DeclareMathSymbol{\qof}{\mathord}{hebrewletters}{114}
\DeclareMathSymbol{\shin}{\mathord}{hebrewletters}{152}

\begin{document}
\author{Peter B. Shalen}
\address{Department of Mathematics, Statistics, and Computer Science (M/C 249)\\  University of Illinois at Chicago\\
  851 S. Morgan St.\\
  Chicago, IL 60607-7045} \email{shalen@math.uic.edu}

\title
[Euler characteristics, lengths of loops, and Wilson's Freiheitssatz]
{
Euler characteristics,  lengths of loops in
hyperbolic $3$-manifolds, and Wilson's Freiheitssatz
}

\begin{abstract}
Let $p$ be a point of  an orientable hyperbolic $3$-manifold
$M$, and let $m\ge1$ and $k\ge2$ be
integers. Suppose that $\alpha_1,\ldots,\alpha_m$ are loops based at
$p$ having length less than $\log(2k-1)$. We show
that  if $G$ denotes the subgroup
of $\pi_1(M,p)$ generated by $[\alpha_1],\ldots,[\alpha_m]$, then 
$\chibar(G)\doteq-\chi(G)\le k-2$; here $\chi(G)$ denotes the Euler
characteristic of the group $G$, which is always defined in this
situation.

This result is deduced from a  result about an arbitrary
finitely generated subgroup $G$ of the fundamental group of an
orientable hyperbolic $3$-manifold. If $\Delta$ is a finite generating
set for $G$, we define the {\it index of freedom} $\iof(\Delta)$ to be
the largest integer $k$ such that $\Delta$ contains $k$ elements that
freely generate a rank-$k$ free subgroup of $G$. We define the {\it
  minimum index of freedom} $\miof(G)$ to be $\min_{\Delta }\iof(\Delta )$, where $\Delta $ ranges over all
finite generating sets for $G$. The result is that
$\chibar(G)<\miof(G)$. The author has recently learned that this is
equivalent to a
special case of a theorem about arbitrary groups due to J. S. Wilson.

\end{abstract}

\maketitle

\section{
Preface to the second version
}

After posting the first version of this paper, which was entitled ``Euler characteristics, free subgroups, and lengths of loops in
hyperbolic $3$-manifolds,'' I was informed by Jack Button that Theorem
B, the proof of which occupies the bulk of the paper, is essentially a special
case of a theorem due to J. S. Wilson. 

Wilson's theorem, which was first proved in \cite{wilson-pro-p},
asserts that if $G$ is a group defined by a presentation with $n$ generators and $m$
relations, where $n>m$, then any generating set for $G$ contains
$n-m$ elements that freely generate a free group of rank $n-m$. In the
language of the present paper, the conclusion says that the minimal index of
freedom $\miof(G)$ is at least $n-m$. 

The {\it deficiency} of a group presentation is defined to be $n-m$,
where $n$ is the number of generators in the presentation and $m$ is
the number of relations. The {\it deficiency} of a finitely
presentable group $G$, which I will denote by $\defish(G)$, is defined to be the maximum of the
deficiencies of all presentations of $G$. (The maximum exists
because the deficiency of any presentation of $G$ is bounded above
by the number of infinite cyclic summands of the commutator quotient
of $G$.) Thus Wilson's result may be stated as follows:

\begin{wilfrei}
For any finitely generated group $G$ we have $\miof(G)\ge\defish(G)$.
\end{wilfrei}

If $G$ is a finitely generated subgroup of the fundamental group of an orientable
hyperbolic $3$-manifold $M$, it follows from Proposition \ref{cocorico}
below, applied to the covering space defined by $M$, that $G$ is
isomorphic to the fundamental group of a compact, orientable,
aspherical $3$-manifold $N$ (so that in particular $G$ is finitely
presented and the Euler characteristic $\chi(G)$ is defined,
cf. Subsection \ref{every little breeze} below).
It follows from a well known result due to D. B. A. Epstein
\cite{epstein}, applied to the compact, orientable, aspherical manifold $N$, that $\defish(G)=1-\chi(G)$,
i.e. $\chibar(G)=\defish(G)-1$ in the notation of this paper. Combining
the latter equality with Wilson's Freiheitssatz immediately gives
Theorem B of this paper.

While Wilson's orginal proof of his Freiheitssatz used the theory of
pro-$p$ groups, he gave another proof in \cite{wilson-elem}
which is quite elementary. It is also remarkably beautiful, involving
surprising arguments about orderable groups, embeddings of group rings
in skew fields, and the so-called Magnus embedding. Although I highly
recommend the paper \cite{wilson-elem}, it is possible that some
readers will be interested in seeing a proof of Theorem B using
$3$-manifold theory and Kleinian groups, which matches the context in
which Theorem B
is applied to the proof of Theorem A. I have therefore decided to
leave the paper on the arXiv, in its present revised form. The present
version differs from
the previous version only in having a new title, this preface, and an
expanded bibliography and abstract.
At some point I may post a new version in which the content of this
preface is integrated into the body of the paper.

%\Pi

\section{
Introduction
}

\Definitions\label{type F def}
We shall say that a space $X$ is  {\it homologically finite} if the
abelian group $\bigoplus_{i=1}^\infty H_i(X;\ZZ)$ is finitely
generated. A homologically finite space has a well-defined
Euler characteristic, which we will denote  by $\chi(X)$. We will
write $\chibar(X)=-\chi(X)$.

A group
$\Pi$ is said to be {\it homologically finite} if $K\doteq K(\Pi,1)$
is homologically finite. In this case we write  $\chi(\Pi)=\chi(K)$ and
$\chibar(\Pi)=\chibar(K)=-\chi(\Pi)$.
\EndDefinitions

It is a well-known fact, which will be  explained in \ref{every
  little breeze} below, that if $G$ is a finitely generated subgroup of the fundamental group of
a hyperbolic $3$-manifold, then $G$ is homologically finite,
so that $\chibar(G)$ is defined. 

One of the main theorems of this paper is:

\begin{theoremA}\label{geometric}
Let $M$ be an orientable hyperbolic $3$-manifold,
let $p$ be a point of $M$, and let $m\ge1$ and $k\ge2$ be
integers. Suppose that $\alpha_1,\ldots,\alpha_m$ are loops based at
$p$ having length less than $\log(2k-1)$. Let $G$ denote the subgroup
of $\pi_1(M,p)$ generated by $[\alpha_1],\ldots,[\alpha_m]$ (so that
$\chibar(G)$ is defined by \ref{every little breeze}). Then
$\chibar(G)\le k-2$. 
\end{theoremA}

This theorem will be used in a forthcoming paper with Rosemary Guzman
to obtain a result that gives a
  new bound on the ratio of the dimension of the mod $2$ homology of a
  closed, orientable $3$-manifold to the volume; this result is
  strictly stronger than the one established in \cite{ratioI}, and is
  essentially stronger than the one established in \cite{ratioII}.

The proof of Theorem A depends on our other main result, Theorem below. Theorem B gives a purely
group-theoretical property of Kleinian groups; the statement of Theorem B involves a few other concepts which
we now define.

\Definitions\label{let freedom ring} Let $G$ be a finitely generated group. We will say that
elements $x_1, \ldots, x_k$ of $G$ are {\it independent}
if $x_1,\ldots,x_k$ freely generate a free subgroup of $G$. If $\Delta $
is a finite generating set for $G$, I will define the {\it index of
  freedom} of $\Delta $, denoted $\iof(\Delta )$, to be the largest integer $k$ such that $\Delta $ contains
$k$ independent elements. We define the {\it minimal index of
  freedom} of $G$, denoted $\miof(G)$, by
$$\miof(G)=\min\nolimits_{\Delta }\iof(\Delta ),$$ where $\Delta $ ranges over all
finite generating sets for $G$.
\EndDefinitions

\begin{theoremB}\label{locke's bagel} 
Let $M$ be an orientable hyperbolic $3$-manifold, and let $G$
be a finitely generated subgroup of $\pi_1(M)$ (so that $\chibar(G)$
is defined by \ref{every little breeze}).
Then $\chibar(G)<\miof(G)$.
\end{theoremB}

Under the hypothesis of Theorem B, the group $G$ is isomorphic to a
Kleinian group, and the proof of the theorem immediately reduces to a
statement about a Kleinian group $\Gamma$, which is Proposition \ref{agol bagol} below.
Using a Dehn surgery trick, one can reduce the proof of Proposition \ref{agol bagol}  to
the case in which $\Gamma$ has no rank-$2$ cusp subgroups.
The essential step in the proof of this case
of Proposition \ref{agol bagol}  is to show that if
$x_1,\ldots,x_m$ are generators for the group $\Gamma$,  if we set
$\Pi_r=\langle x_1,\ldots,x_r\rangle\le \Gamma$ for $r=1,\ldots,m$, and if
for a given $r\in\{1,\ldots,m-1\}$ we have
$\chibar(\Pi_{r+1})>\chibar(\Pi_r)$, then $x_{r+1}$ has infinite order
in $\Gamma$, and $\Pi_{r+1}$ is the free
product of its subgroup $\Pi_{r}$ with its infinite cyclic subgroup
$\langle x_{r+1}\rangle$. This is proved by relating the quantities
$\chibar(\Pi_{r})$ for $1\le r\le m$ to the dimensions of certain
varieties of representations. 

For each $r$ the representations of
$\Pi_r$ in $\pizzle(\CC)$ are identified with the points of a complex
affine algebraic set, which for the purposes of this sketch will be
denoted $\barR_r$. The inclusion homomorphism $\Pi_r\to\pizzle(\CC)$,
regarded as a representation, lies in a unique irreducible component
$V_r$ of $\barR_r$, and has complex dimension $3\chibar(\Pi_r)+3$. If
for a given $r<m$ the group
$\Pi_{r+1}$ is {\it not} the free
product of  $\Pi_{r}$ with
$\langle x_{r+1}\rangle$ (or if $x_{r+1}$ has finite order), it can be shown that $V_{r+1}$ is isomorphic to
a proper subvariety of $V_r\times\pizzle(\CC)$; this implies that
$\dim V_{r+1}<3+\dim V_{r}$, and hence that
$\chibar(\Pi_{r+1})\le\chibar(\Pi_{r})$. This completes the sketch of
the proof of Theorem B. (Parts of the argument sketched here are
embodied in the statement and proof of Lemma \ref{people who can't
  count}.) 

The idea of using dimensions of representation varieties in this argument
was suggested by Ian Agol's observation that one can use dimensions of
representation varieties to give alternative proofs of \cite[Theorem
VI.4.1]{JS}, and its generalizations \cite[Appendix, Theorem
A]{threefree} and \cite[Theorem 7.1]{accs}, using dimensions of representation varieties.

Theorem A is proved by combining Theorem B with the so-called
$\log(2k-1)$ Theorem. The latter result, which gives information about
the lengths of loops $\alpha_1,\cdots,\alpha_k$ based at a point $p$
of an orientable hyperbolic $3$-manifold $M$ under the assumption that
the elements $[\alpha_1],\ldots,[\alpha_k]$ freely generate a free
subgroup of $\pi_1(M,p)$, was proved under an additional hypothesis in
\cite{accs}. The general version of the $\log(2k-1)$ Theorem is proved
by combining the results of \cite{accs} with the tameness theorem
\cite{agol}, \cite{cg} and the density theorem \cite{ohshika},
\cite{namazi-souto} for Kleinian groups. The general
version of the $\log(2k-1)$ Theorem is presented as Theorem 4.1 of
\cite{acs-surgery}, which we will quote in the proof of Theorem A.

The proofs of Theorems B and A will be given in Section \ref{da
  proofs}, after needed background about representations,
representaion varieties, character varieties, and deformation spaces
of Kleinian groups, which occupies Section \ref{prelims}.

We will use the following standard conventions. A $3$-manifold $M$ is
said to be {\it irreducible} if $M$ is connected and every (tame)
$2$-sphere in $M$ bounds a $3$-ball. To write $Y\le X$, where $X$ is a
group, means that $Y$ is a subgroup of $X$. If $A$ is a subset of a
group, $\langle A\rangle$ denotes the subgroup generated by $A$.

I am grateful to Steve Boyer, Ken Bromberg and Dick Canary for their
assistance. Boyer, with tremendous patience, helped me navigate the material that is summarized
in Subsection \ref{introducing X bar}. Bromberg and Canary helped me
locate the result from \cite{fps} which is used in the proof of Lemma
\ref{before doodle poodle}. Canary pointed out Theorem 8.44 of
\cite{mishabook}, which is used in the proof of Lemma \ref{sorta marden}.

\section{
Preliminaries
}\label{prelims}

\Definition\label{supreme court}
A {\it compact core} of a $3$-manifold $M$ is a compact
three-dimensional submanifold
$N$ of $M$ such that the inclusion $N\to M$ is a homotopy equivalence.
\EndDefinition

Although the following result is well known, we include a proof for
clarity.

\Proposition\label{cocorico}
Every orientable hyperbolic $3$-manifold with finitely generated fundamental
group has a compact core.
\EndProposition

\Proof
If the given manifold $M$ is isometric to $\HH^3$, any ball in $M$ is
a compact core. If $M$ is not isometric to $\HH^3$,  we observe that
since $M$ admits $\HH^3$ as a non-trivial covering space, $M$ is
irreducible and has infinite fundamental group. We then apply
\cite[Proposition 3.8]{small}, a consequence of the main result of
\cite{core}, which asserts that if $M$ is an irreducible, orientable
3-manifold with finitely generated fundamental group, then
there is a compact, irreducible submanifold $N$ of $M$ such that the
inclusion homomorphism $\pi_1(N)\to\pi_1(M)$ is an isomorphism. (The
hypothesis of finite generation was unfortunately left out of the
statement of \cite[Proposition 3.8]{small}, but the proof given there
establishes the statement given here.) In
particular $\pi_1(N)$ is infinite; this, together with the
irreducibility of $N$, implies that $N$ is aspherical (see
\cite[Proposition 3.6]{small}). As the hyperbolic manifold $M$ is also
aspherical, it now follows that the inclusion $N\to M$ is a homotopy
equivalence.
\EndProof

\Number\label{every little breeze}
Let $G$ be a finitely generated
subgroup of the fundamental group of an orientable hyperbolic
$3$-manifold $M$. The covering space $\tM$ of $M$ corresponding to the
subgroup $G$ has a compact core $N$ by Proposition \ref{cocorico}; in
particular $N$ is aspherical and $\pi_1(N)\cong G$. Since $N$ is
compact, it follows that $G$ is homologically finite.

This establishes the fact, which was mentioned in the introduction and
is necessary background for Theorems A and B, that a finitely generated
subgroup $G$ of the fundamental group of an orientable hyperbolic
$3$-manifold is homologically finite, so that $\chibar(G)$ is defined.
\EndNumber

\Lemma\label{before doodle poodle}
If $M$ is an orientable hyperbolic $3$-manifold with finitely
generated Kleinian group, there exist an orientable hyperbolic
$3$-manifold $M_0$ with no rank-$2$ cusps such that
$\chibar(M_0)=\chibar(M)$, and a surjective homomorphism
from $\pi_1(M)$ to $\pi_1(M_0)$. (Here $\chibar(M)$ and
$\chibar(M_0)$ are defined in view of \ref{every little breeze}.)
\EndLemma

\Proof
It follows from Proposition \ref{cocorico} that $\dim
H_2(M;\QQ)<\infty$. The number of ends of the manifold $M$ is at most $1+\dim
H_2(M;\QQ)$, and is in particular finite. We shall prove the result by
induction on the number of ends of $M$.

If the number of ends of $M$ is $0$, then $M$ is compact and therefore
has no cusps. The assertion is trivial in this case, since we can take
$M_0=M$ and choose the identity map as required homomorphism. Now
suppose that the number of ends of $M$ is $m>0$, and that the assertion is
true for manifolds with $m-1$ ends. If $M$ has no rank-$2$ cusps, the
assertion is again trivial. Now suppose that $M$ has at least one
rank-$2$ cusp, and fix a submanifold $H$ of $M$ which is a standard
neighborhood of a rank-$2$ cusp.

The $2$-manifold $\partial H$ is a torus, and inherits a Euclidean
metric from the hyperbolic metric on $M$. For each homotopically
non-trivial simple closed curve
$s$ in $\partial H$, we denote by $Z(s)$ the orientable $3$-manifold, defined up
to homeomorphism and dependent only on the isotopy class of $s$, which is obtained
from the disjoint union $\overline{M-H}\discup (D^2\times S^1)$ by
gluing $\partial(\overline{M-H})=\partial H$ to $\partial (D^2\times
S^1)=(\partial D^2)\times S^1$ by a homeomorphism which maps $s$ to
$(\partial D^2)\times\{\text{pt}\}$. (For any $s$, the manifold
$Z(s)$ is said to be obtained from $\overline{M-H}$ by Dehn filling.)

For any closed geodesic $s$ in $\partial H$, since $\overline{M-H}$
is homotopy equivalent to $M$, we have
$\chibar(Z(s)=\chibar(\overline{M-H})-\chibar(\partial
H)+\chibar(D^2\times S^1)=\chibar(M)-0+0=\chibar(M)$. Furthermore, the
inclusion homomorphism $\pi_1(\overline{M-H})\to\pi_1(Z(s))$ is
surjective, and hence there is a surjective homomorphism from
$\pi_1(M)$ to $\pi_1(Z(s))$.

We now apply the case $k=1$ of
\cite[Theorem 1.4]{fps}, which asserts that if the closed geodesic in
$\partial H$ representing the isotopy class of $s$ has length strictly
greater than $6$, then $Z(s)$ admits a hyperbolic metric. Since there
are only finitely many isotopy classes of geodesics in $\partial H$
whose lengths are bounded above by a given constant, we may fix a
simple closed curve $s_1$ in $\partial H$ such that $Z(s_1)$ is
homeomorphic to a hyperbolic $3$-manifold. 

The manifold $Z(s_1)$ has $m-1$ ends. Hence, by the induction
hypothesis, there exist  an orientable hyperbolic
$3$-manifold $M_0$ with no rank-$2$ cusps such that
$\chibar(M_0)=\chibar(Z(s_1))$, and a surjective homomorphism
from $\pi_1(Z(s_1))$ to $\pi_1(M_0)$. Since
$\chibar(Z(s_1))=\chibar(M)$, we have
$\chibar(M_0)=\chibar(M)$. Furthermore, since there is a surjective homomorphism
from $\pi_1(M)$ to $\pi_1(Z(s_1))$, we obtain by composition a
surjective homomorphism from
$\pi_1(M)$ to $\pi_1(M_0)$. This completes the induction.
\EndProof

\Lemma\label{doodle poodle}
If $\Gamma$ is a finitely generated Kleinian group, there exist
a finitely generated Kleinian
group $\Gamma_0$ with no rank-$2$ cusp subgroups, with
$\chibar(\Gamma_0)=\chibar(\Gamma)$, and a surjective homomorphism
$\eta:\Gamma\to\Gamma_0$. (Here $\chibar(\Gamma)$ and
$\chibar(\Gamma_0)$ are defined by \ref{kleinian def}.)
\EndLemma

\Proof
We apply Lemma \ref{before doodle poodle} to the orientable hyperbolic
$3$-manifold $M\doteq\HH^3/\Gamma$. The latter lemma gives an  orientable hyperbolic
$3$-manifold $M_0$ with no rank-$2$ cusps such that
$\chibar(M_0)=\chibar(M)$, and a surjective homomorphism
from $\pi_1(M)$ to $\pi_1(M_0)$. We may write $M_0=\HH^3/\Gamma_0$
where $\Gamma$ is a Kleinian group without rank-$2$ cusp
subgroups. We have $\Gamma\cong\pi_1(M)$ and
$\Gamma_0\cong\pi_1(M_0)$; and since the hyperbolic $3$-manifolds $M$
and $M_0$ are aspherical, we have
$\chibar(\Gamma_0)=\chibar(M_0)=\chibar(M)=\chibar(\Gamma)$. The
assertion follows.
\EndProof

The next two subsections involve the notion of index of freedom and
minimal index of freedom, which were defined in the introduction to
this paper.

\Number\label{well duh}Note that if a finitely generated group $G$ is
non-trivial and torsion-free, then any generating set for $G$ contains
an element of infinite order, and hence $\miof(G)\ge1$.        
\EndNumber

\Lemma\label{quotient miof}
Let $G_1$ and $G_2$ be groups, and suppose that there is a
surjective homomorphism $\eta:G_1\to G_2$. Then
$\miof(G_2)\le\miof(G_1)$. Furthermore, for any generating set $\Delta$ for
$G_1$, we have $\iof(\eta(\Delta))\le\iof(\Delta)$.
\EndLemma

\Proof
We prove the second assertion first. Write
$\Delta=\{x_1,\ldots,x_m\}$. Set $k=\iof(\eta(\Delta))$. Then after
reindexing the $x_i$ we may assume that $\eta(x_1),\ldots,\eta(x_k)$
are independent. It follows that $x_1,\ldots,x_k$ are independent, and
hence that $\iof(\Delta)\ge k$, as required. To prove the first
assertion, note that if $\Delta$ is an arbitrary generating set for
$G_1$, the first assertion gives
$\iof(\Delta)\ge\iof(\eta(\Delta)\ge\miof(G_2)$; hence 
$\miof(G_1)\ge\miof(G_2)$.
\EndProof

%\Lemma \label{listen mister bill be}Let $\Gamma_1$ and $\Gamma_2$ be finitely generated subgroups of
%$\zzle(\CC),$ neither of which contains any non-trivial scalar matrix. Then there is an element $T$ of
%$\ggle(\CC)$ such that the subgroup generated by $\Gamma_1$ and
%$T\Gamma_2T^{-1}$ is a free product of its subgroups $\Gamma_1$ and $T\Gamma_2T^{-1}$.
%\EndLemma

%\Proof
%Since the $\Gamma_i$ are finitely
%generated, there is a finitely generated
%extension $F\subset\CC$ of $\QQ$ such that $\Gamma_1$ and $\Gamma_2$ are
%contained in $\zzle(F)$. The proof of Theorem 3 of \cite{wehrfritz}
%(see also \cite{nisnevic})
%then shows that if $X,Y,Z,W$ are independent transcendentals over $F$,
%and if we set 
%$$T_0=\begin{pmatrix}X&Y\\Z&W\end{pmatrix},$$
%then  the subgroup of $\zzle(F(X,Y,Z,W))$ generated by $\Gamma_1$ and
%$T_0\Gamma_2T_0^{-1}$ is a free product of its subgroups $\Gamma_1$
%and $T_0\Gamma_2T_0^{-1}$. Since $F$ is finitely generated, the
%inclusion of $F$ in $\CC$ extends to an embedding of the field
%$F(X,Y,Z,W)$ in $\CC$, and the conclusion follows.  
%\EndProof

The following lemma is a very special case of a result that is implicit in \cite{nisnevic}. The latter result has been refined in  \cite{wehrfritz} and elsewhere, and was rediscovered in \cite{shalen-amalgam}. As an explicit statement of the lemma is not easy to find, we have provided a simple, self-contained proof.

\Lemma \label{no but it will be}Let $\Gamma$ be a finitely generated subgroup of
$\ggle(\CC)$ containing no non-trivial scalar matrix. 
Then there is an element $L$ of
infinite order in $\zzle(\CC)$ such that the subgroup  $\langle
L,\Gamma\rangle$ of $\ggle(\CC))$ is a free product of $\Gamma$ with the infinite cyclic
group $\langle L\rangle$, and contains no non-trivial scalar matrix.
\EndLemma

\Proof
Since $\Gamma$ contains no non-trivial scalar matrix, each non-trivial
element of $\Gamma$ has one or two eigenvectors in $\CC^2$. Since
$\Gamma$ is countable and $\CC^2$ is not, there is a non-zero vector in $\CC^2$ which is not an eigenvector of any element of $\Gamma$. Hence after modifying $\Gamma$ by a conjugation, we may assume that $\Gamma$ contains no upper triangular matrix except the identity.

Let $F$ denote the countable subfield of $\CC$ generated by the entries of the elements of $\Gamma$, and fix an element $t$ of $\CC$ which is transcendental over $F$.
We shall show that the conclusion of the lemma holds if we set $L=\begin{pmatrix}1&t\\0&1\end{pmatrix}$. This is equivalent to the assertion that if $F(X)$ denotes the rational function field in one indeterminate over $X$, and if we set $\Lambda=\begin{pmatrix}1&X\\0&1\end{pmatrix}$, then the subgroup of $\ggle(F[X])\le\ggle(F(X))$ generated by 
$\Gamma$ and $\Lambda$ is a free product of $\Gamma$ with the infinite cyclic
group $\langle \Lambda\rangle$, and contains no non-trivial scalar matrix.

We consider the abstract free product $\Gamma\star\langle \Lambda\rangle$, and the homomorphism $h: \Gamma\star\langle \Lambda\rangle\to
\ggle(F[X])$ which is restricts to the inclusion on each factor; we
must show that $h$ does not map any non-trivial element of
$\Gamma\star\langle \Lambda\rangle$ to a scalar matrix. Any non-trivial element of $\Gamma\star\langle \Lambda\rangle$ is conjugate either to an element of a factor or to an element of the form
\Equation\label{tulip}
w=\gamma_1\Lambda^{m_1}\cdot \gamma_k\Lambda^{m_k},
\EndEquation
where $k\ge1$, the elements $\gamma_1,\ldots,\gamma_k$ of $\Gamma$ are
non-trivial, and $m_1,\ldots,m_k$ are non-zero integers. 
By hypothesis $\Gamma$ contains no non-trivial scalar matrix, and it
is clear from the definition of $\Lambda$ that $\langle \Lambda\rangle$ also
contains no non-trivial scalar matrix. It therefore suffices to prove
that an element $w$ of the form (\ref{tulip}), where $k\ge1$ and the
$\gamma_i$ and $m_i$ satisfy the conditions stated above, cannot be
mapped by $h$ onto a scalar matrix. 
By induction on $k\ge1$ we shall establish the following stronger assertion:
\Claim\label{toodly-froodly}
If $k$ is a strictly positive integer, and $w\in \Gamma\star\langle \Lambda\rangle$ is given by \ref{tulip}, where
$k\ge1$, the elements $\gamma_1,\ldots,\gamma_k$ of $\Gamma$ are
non-trivial, and $m_1,\ldots,m_k$ are non-zero integers, then $h(w)$
has the form $\begin{pmatrix}A&B\\C&D\end{pmatrix}\in\ggle(F[X])$,
where $A$ and $C$ are polynomials of degree at most $k-1$, while $B$
is a polynomial of degree at most $k$ and $D$ is a polynomial of
degree exactly $k$. (In particular we have $A\ne D$, so that $h(w)$
cannot be a scalar matrix.)
\EndClaim

In the statement of \ref{toodly-froodly}, the polynomial $0$ is understood to have degree $-\infty$.

In the case $k=1$ of \ref{toodly-froodly}, we have $w=\gamma \Lambda^{m }$, where $m $ is a non-zero integer and $\gamma =
\begin{pmatrix}a &b \\c &d \end{pmatrix}$ is a non-trivial element of $\Gamma$. This gives
$$h(w)=\begin{pmatrix}a &b \\c &d \end{pmatrix}
\begin{pmatrix}1&m X\\0&1\end{pmatrix}
=\begin{pmatrix}a &a m X+b \\c &c m X+d \end{pmatrix}.
$$
Thus the entries in the first column of $h(w)$ are elements of $F$,
while the upper right-hand entry is a polynomial of degree at most
$1$. Furthermore, since $\Gamma$ contains no upper triangular matrix
except the identity, we have $c\ne0$, and hence the lower right-hand entry of $h(w)$ has degree exactly $1$. This establishes the base case for the inductive proof of \ref{toodly-froodly}.

For the induction step, suppose that $w=\gamma_1\Lambda^{m_1}\cdots \gamma_k\Lambda^{m_k}$ satisfies the hypothesis of \ref{toodly-froodly} for a given $k>1$, set $w^*=\gamma_1\Lambda^{m_1}\cdots \gamma_{k-1}\Lambda^{m_{k-1}}$, and assume  that  $h(w^*)=\begin{pmatrix}A^*&B^*\\C^*&D^*\end{pmatrix}$, where $A^*$ and $C^*$ have degree at most $k-2$, while $B^*$ has degree at most $k-1$, and $D^*$ has degree exactly $k-1$. 
As the base case has been established, we may write
$h(
\gamma_k\Lambda^{m_k}
)=\begin{pmatrix}A^\dagger&B^\dagger\\C^\dagger&D^\dagger\end{pmatrix}$, where $A^\dagger$ and $C^\dagger$ have degree at most $0$, while $B^\dagger$ has degree at most $1$, and $D^\dagger$ has degree exactly $1$. Then
$$h(w)=h(w^*) h(
\gamma_k\Lambda^{m_k})=\begin{pmatrix}A^* A^\dagger+B^*C^\dagger&
A^* B^\dagger+B^*D^\dagger\\
C^* A^\dagger+D^*C^\dagger&
C^* B^\dagger+D^*D^\dagger\end{pmatrix}.
$$
%\star
It now follows that the entries in the first column of $h(w)$ have degree at most $k-1$, while its upper right-hand entry has degree at most $k$. In the lower right-hand entry  $C^* B^\dagger+D^*D^\dagger$, the first term has degree at most $k-1$ and the second term has degree exactly $k$; hence this lower right-hand entry has degree exactly $k$, and the induction is complete.
%Since $\Gamma$ has trivial center, it contains no scalar matrices
%except the identity. Let $A_0$ be any element of $\zzle(\CC)$ which has
%infinite order and is not a scalar matrix. According to Lemma
%\ref{listen mister bill be}, applied with $\Gamma_1=\Gamma$ and
%$\Gamma_2=\langle A_0\rangle$, there is an element $T$ of $\ggle(\CC)$
%such that the subgroup generated by $\Gamma$ and
%$\langle TA_0T^{-1}\rangle$ is a free product of its subgroups $\Gamma$ and
%$\langle TA_0T^{-1}\rangle$. The result follows upon setting $A=TA_0T^{-1}$.
\EndProof

\Number\label{introducing X bar}
The rest of this section will involve the use of representation
varieties. We will  take the point of view used in \cite{bznorms}.

By a {\it representation} $\barrho$ of a group $\Gamma$ in
$\pizzle(\CC)$, we mean simply a homomorphism from $\Gamma$ to
$\pizzle(\CC)$. A representation of $\Gamma$ is
said to be {\it reducible} if there is a $1$-dimensional subspace of
$\CC^2$ which is invariant under $\barrho(\Gamma)$. Otherwise it is said
to be {\it irreducible}.

 In \cite[Section
3]{bznorms}, it is shown how to identify $\pizzle(\CC)$ with an
algebraic set in some complex affine space; and it is pointed out that, if $\Gamma$ is a group with a given
  finite generating set $S$, and if one identifies an arbitrary
  representation $\barrho$ of $\Gamma$ in $\pizzle(\CC)$ with the point
  $(\barrho(x))_{x\in S}\in \pizzle(\CC)^S$, then the set of all
  $\pizzle(\CC)$-representations of  $\Gamma$ is identified with a
  complex affine
  algebraic subset $\barR(\Gamma)$ of $\pizzle(\CC)^S$.
% denotes a complex affine
%  algebraic whose points are identified
%  with the representations of $\Gamma$ in $\zzle(\CC)$. Concretely, we may regard $\barR(\Gamma)$
 % as an algebraic subset of $\pizzle(\CC)^S$, where; and the point of
%$\barR(\Gamma)$ identified with a given representation $\barrho$ is
%$(\barrho(x))_{x\in S$.

As in \cite{bznorms}, we will also consider the
$\pizzle(\CC)$-character variety $\barX(\Gamma)$ for an arbitrary
finitely generated group $\Gamma$. This is the analogue of the
$\zzle(\CC)$-character variety of $\Gamma$ that was considered in
\cite {splittings}, and we will review its definition and some of its useful
properties here. (We will need to use $\barX(\Gamma)$ in
this paper because Theorem 8.44 of \cite{mishabook}, which is quoted
in the proof of Lemma \ref{sorta marden}, is stated in terms of
$\barX(\Gamma)$.)

There is a natural action $(g,\barrho)\mapsto \g\cdot\barrho$ of $\pizzle(\CC)$ on
$\barR(\Gamma)$, where the representation $g\cdot\barrho$ is defined by
setting $g\cdot\barrho(x)=g\barrho(x)g^{-1}$ for each $x\in\Gamma$. This
action will be referred to as the action by {\it conjugation}, and its
orbits  will be called {\it conjugacy classes} of
representations. In \cite[Section 3]{bznorms},  $\barX(\Gamma)$ is
defined as the quotient, in the category of complex affine algebraic
sets, of $\barR(\Gamma)$ by this action. There is a surjective morphism of affine
algebraic sets $\bart: \barR(\Gamma)\to \barX(\Gamma)$ which is
constant on each conjugacy class of representations. For each
$\barrho\in\barR(\Gamma)$, we call $\bart(\barrho)$ the {\it character} of $\barrho$.

(As in \cite{bznorms}, the use of bars in $\barR(\Gamma)$,
$\barX(\Gamma)$ and $\bart$, and the use of $\barrho$ as the default
notation for an element of $\barR(\Gamma)$, are meant to emphasize
that we are dealing with representations in $\pizzle(\CC)$ rather than
$\zzle(\CC)$.
However, we have avoided using $\chibar_{\barrho}$ to mean
$\bart(\barrho)$, as is done in \cite{bznorms}, since the symbol
$\chibar$ figures prominently in the present paper with a different meaning.)

In the discussion at the beginning of Section 3 of \cite{bznorms}, the
authors quote 
%the general theory presented in 
\cite{newstead} and \cite{jm} 
%involving actions of geometrically reductive algebraic groups on
%algebraic sets, 
for the precise definition and the basic properties of
the quotient object $\barX(\Gamma)$. The theory presented in
\cite{newstead} applies to an arbitrary action (in the category of
algebraic sets) of $\pizzle(\CC)$ on an arbitrary algebraic set,
whereas in \cite{jm} the emphasis is on the specific case of the
action on $\barR(\Gamma)$ by conjugation, where $\Gamma$ is a finitely
generated group. (In both of these sources, the role of $\pizzle(\CC)$
is played by a more general algebraic group, but we will implicitly
specialize to the case of $\pizzle(\CC)$ in the following discussion.)

According to the definition given on page 53 of \cite{jm}, a point
$\barrho$ of $\barR(\Gamma)$ is {\it stable} if and only if its
conjugacy class is
Zariski-closed and its stabilizer (isotopy subgroup) under the action
of $\pizzle(\CC)$ by conjugation is
finite.

We will need the following properties of $\barX(\Gamma)$ and $\bart$:

\Claim\label{zeroth base}
A point of $\barrho\in\barR(\Gamma)$ is
stable if and only if $\barrho$ is an
irreducible representation. 
\EndClaim

\Claim\label{first base}
For each irreducible representation $\barrho\in\barR(\Gamma)$, the set
$\bart^{-1}(\bart(\barrho))$ is precisely the conjugacy class of the
representation $\barrho$.
\EndClaim

\Claim\label{second base}
Every  Zariski-closed subset
 of $\barR(\Gamma)$ which is invariant under the action of
 $\pizzle(\CC)$ by conjugation is mapped by $\bart$ onto  a Zariski-closed subset of
 $\barX(\Gamma)$. 
\EndClaim

\Claim\label{third base}
The set of all characters of irreducible representations of $\Gamma$
in $\pizzle(\CC)$ is a Zariski-open subset of $\barX(\Gamma)$.
\EndClaim

Assertion \ref{zeroth base} follows from \cite[Theorem 1.1]{jm}.
On p. 753 of \cite{bznorms}, it is pointed out that Assertion
\ref{first base} follows from \cite[Corollary 3.5.2]{newstead} upon
combining it with Assertion \ref{zeroth base}. In the proof
of \cite[Lemma 4.1]{bznorms}, it is pointed out that Assertion
\ref{second base} follows from \cite[Theorem
3.3.5(iv)]{newstead}. 

To prove \ref{third base}, we first note that in view of \ref{zeroth
  base} it is enough to prove that $\bart$ maps the set $S$ of stable
points of $\barR(\Gamma)$ onto a Zariski-open subset of $\barX(\Gamma)$. It
is pointed out on
page 54 of \cite{jm}, where it is deduced from Proposition 3.8 of
\cite{newstead}, that $S$ is Zariski-open in $\barR(\Gamma)$. Hence
$S'=\barR(\Gamma)-S$ is Zariski-closed. Clearly $S$ and $S'$ are invariant under the
action of $\pizzle(\CC)$. It follows from Assertion \ref{first base}
that $\bart(S)\cap\bart(S')=\emptyset$; since $\bart$ is surjective,
it now follows that $\bart(S)=\barR(\Gamma)-\bart(S')$. Since $S'$ is
Zariski-closed in $\barR(\Gamma)$, Assertion \ref{second base}
implies that $\bart(S')$ is Zariski-closed in $\barX(\Gamma)$; hence
$\bart(S)$ is Zariski-closed in $\barX(\Gamma)$, as required.

\EndNumber

\Lemma\label{before sorta marden}
Let $\Gamma$ be a finitely generated group, and let
$\barrho_0\in\barR(\Gamma)$ be an irreducible representation of $\Gamma$
in $\pizzle(\CC)$. Suppose that $\bart(\barrho_0)$ lies in a unique
component $\barX_0$ of
$\barX(\Gamma)$, and set $d=\dim\barX_0$. Then $\barrho_0$
lies in a unique irreducible component $\barR_0$ of $\barR(\Gamma)$, and $\dim\barR_0=d+3$.
\EndLemma

\Proof[Proof (cf. proof of Lemma 4.1 in \cite{bznorms})]
Since $\barX_0$ is an irreducible component of
$\barX(\Gamma)=\bart(\barR(\Gamma))$,
%According to 
%\redcomment{
%reference,
%}
there is an irreducible component $\barR_0$ of $\barR(\Gamma)$ such
that 
%$\bart(\barR_0)=\barX_0$. 
$\bart(\barR_0)$ is a Zariski-dense subset of $\barX_0$. 
If we define a map
$h:\barR_0\times\pizzle(\CC)\to \barR(\Gamma)$ by $h(\barrho,g)$ to be
the representation $x\to g\barrho(x)g^{-1}$, 
then since
$\barR_0\times\pizzle(\CC)$ is irreducible,
$h(\barR_0\times\pizzle(\CC))$ must be contained in a single component
of $\barR(\Gamma)$; but we have $h(\barR_0\times\pizzle(\CC))\supset
h(\barR_0\times\{1\}=\barR_0$, and hence
$h(\barR_0\times\pizzle(\CC))=\barR_0$. This shows that $\barR_0$ is
invariant under the action of $\pizzle(\CC)$ on $\barR(\Gamma)$ by
conjugation. Since the irreducible component $\barR_0$ of
$\barR(\Gamma)$  is also Zariski-closed in  $\barR(\Gamma)$, It follows from \ref{second base} that
$\bart(\barR_0)$ is a Zariski-closed subset of
 $\barX(\Gamma)$. Hence $\bart(\barR_0)=\barX_0$. 

Let $\barX_1$ denote the union of all irreducible components of
$\barX(\Gamma)$ that are distinct from $\barX_0$, and set
$U_1=\barX(\Gamma)-\barX_1$. Then $U_1$ is Zariski-open in
$\barX(\Gamma)$, and is contained in $\barX_0$. On the other hand, if
$U_2\subset \barX(\Gamma)$ denotes the set of all characters of
irreducible representations of $\Gamma$ in $\pizzle(\CC)$, then $U_2$
is Zariski-open in $ \barX(\Gamma)$ by \ref{third base}.
Hence $U\doteq U_1\cap U_2$ is Zariski-open in
$\barX(\Gamma)$, and is contained in $\barX_0$. 

Set $\barx_0=\bart(\barrho_0)$.
Since by hypothesis
$\barrho_0$ is an irreducible representation and $\barX_0$ is the only irreducible
component of $\barX(\Gamma)$ containing $\barx_0$, we have $\barx_0\in U$.

Consider an arbitrary point $\barx\in U$. Since $\barx\in U_2$, it follows
from \ref{first base}
that $\bart^{-1}(\barx)$ is a single conjugacy class of
irreducible representations. Since $\barx\in
U_1\subset\barX_0=\bart(\barR_0)$, 
and
since $\barR_0$ is invariant under conjugation, it follows that
$\bart^{-1}(\barx)\subset\barR_0$. Hence we have
$\bart^{-1}(U)\subset\barR_0$. But $\bart^{-1}(U)$ is Zariski-open in
$\barR(\Gamma)$ since $U$ is Zariski-open in $\barX(\Gamma)$, and since $\bart(\barrho_0)=\barx_0\in U$,
we have $\barrho_0\in\bart^{-1}(U)$. Thus
$\bart^{-1}(U)$ is a Zariski neighborhood of $\barrho_0$ in
$\barR(\Gamma)$, and is contained in $\barR_0$; this shows that
$\barR_0$ is the only irreducible component of
$\barR(\Gamma)$ containing $\barrho_0$. It remains to show that
$\dim\barR_0=d+3$.

We have observed that for each
$\barx\in U$, the fiber $\bart^{-1}(\barx)$ is a single conjugacy class of
irreducible representations. If we fix a representation
$\barrho\in\bart^{-1}(\barx)$, the algebraic set $\bart^{-1}(\barx)$ is
isomorphic to the coset space $\pizzle(\CC)/C$, where $C$ denotes the
stabilizer of $\barrho(\Gamma)$ under the action of $\pizzle(\CC)$ by
conjugation. Since $\barrho$ is a stable point by \ref{zeroth base}, it
follows from the definition of stability (see \ref{introducing X bar})
that
%According to
%Lemma \ref{funny}, 
$C$ is a finite group, and hence
%reference,
%the image of an irreducible representation of $\Gamma$ has trivial
%centralizer in $\pizzle(\CC)$, and hence the action of $\pizzle(\CC)$
%on a fiber $\bart^{-1}(\barx)$ is free. We therefore have 
$\dim
\bart^{-1}(\barx)=\dim\pizzle(\CC)=3$; that is, all fibers of the map
$\bart|\bart^{-1}(U): \bart^{-1}(U)\to U$ are three-dimensional. Hence
$\dim \bart^{-1}(U) =3+\dim U$. But
$U$ is in particular a Zariski-open subset of $\barX_0$,
and is non-empty since it contains $\barx_0$; hence $\dim U=\dim
X_0=d$. and so $\dim \bart^{-1}(U) =d+3$. Finally, since
$\bart^{-1}(U)$ is in particular a Zariski-open subset of the
irreducible variety $\barR_0$,
we have $\dim\barR_0=\dim\bart^{-1}(U)=d+3$.
\EndProof

\DefinitionRemark\label{kleinian def}
By a {\it Kleinian group} we mean a discrete subgroup of
$\pizzle(\CC)$.

If $\Gamma$ is
a finitely generated torsion-free Kleinian group, then $\Gamma$ is
isomorphic to the
fundamental group of the orientable hyperbolic $3$-manifold
$\HH^3/\Gamma$, and is therefore homologically finite by \ref{every little
  breeze}. In particular, $\chibar(\Gamma)$ is defined for every finitely generated torsion-free Kleinian group $\Gamma$.
\EndDefinitionRemark

\Lemma\label{sorta marden}
Let $\Gamma$ be a finitely generated, torsion-free group, and
let $\barrho_0$ be a faithful representation of $\Gamma$ in $\pizzle(\CC)$
such that  $\barrho_0(\Gamma)$ is a Kleinian group
(so that $\chibar(\Gamma)$ is defined by \ref{kleinian def}). 
Assume that $\Gamma$ has no rank-$2$ 
cusp subgroups.
Then the point
$\barrho_0$ of $\barR(\Gamma)$ lies in a unique irreducible component
$\barR_0$ of $\barR(\Gamma)$, and $\barR_0$
has dimension $3\chibar(\Gamma)+3$. 
\EndLemma

\Remark\label{can't be bothere}
It is likely that by using a little more geometric invariant theory, one can
prove that under the hypothesis of Lemma \ref{sorta marden},
the point  $\barrho_0$ of $\barR(\Gamma)$ is smooth. As we do not
need this stronger result, we have not attempted to include a
proof of it.
\EndRemark

\Proof[Proof of Lemma  \ref{sorta marden}]
%Let $M$ denote the hyperbolic $3$-manifold $\HH^3/\Gamma$. Since $M$
%is covered by $\HH^3$, it is irreducible. Since in addition
%$\pi_1(M)\cong\Gamma$ is finitely generated, it follows from 
%\redcomment{
%reference, pointing out that I'm using Kapovich's notation
%}
%that there is a compact three-dimensional submanifold $M_c$ of $M$ such that the inclusion
%$\pi_1(M_c)\to\pi_1(M)$ is a homotopy equivalence. In particular $M_c$ is
%a compact, aspherical manifold, so that $\Gamma\cong\pi_1(M_c)$ is of
%type F. This is the first assertion of the lemma.
%To prove the remaining assertions, 
First consider the case in which
the Kleinian group
$\barrho_0(\Gamma)\cong\Gamma$ is elementary. An elementary, torsion-free Kleinian group
without rank-$2$ cusp subgroups is either trivial or infinite cyclic. If $\Gamma$
is trivial then $\chibar(\Gamma)=-1$ and $\barR(\Gamma)$ is a single
point; if $\Gamma$
is infinite cyclic, then $\chibar(\Gamma)=0$ and the variety
$\barR(\Gamma)$  is isomorphic to $\pizzle(\CC)$, and is therefore an
irreducible complex affine variety of dimension $3$. Thus the
conclusions hold in the elementary case.

Now suppose that $\Gamma$ is non-elementary. According to
\cite[Theorem 8.44]{mishabook}, $\bart(\barrho_0)$ is a smooth point of
$\barX(\Gamma)$, and hence lies in a unique irreducible component
$\barX_0$ of $\barX(\Gamma)$; and furthermore, we have $\dim
X_0=3\chibar(\partial M_c)/2= 3\chibar(M_c)=3\chibar(\Gamma)$, where
$M_c$ denotes a compact core of $M$ (see Definition \ref{supreme
  court} and Proposition \ref{cocorico}).
(The quantity denoted $\tau$ in the statement of \cite[Theorem
8.44]{mishabook} is the number of conjugacy classes of rank-$2$ cusp subgroups of
$\Gamma$, which according to the hypothesis of the present lemma is equal to $0$.)
Furthermore, since $\barrho_0(\Gamma)$ is non-elementary, the
representation $\barrho_0$ is
irreducible. It now follows from Lemma 
 \ref{before sorta marden} that $\barrho_0$ lies in a unique component
 $\barR_0$ of $\barR(\Gamma)$, and that $\dim\barR_0=\dim\barX_0+3=3\chibar(\Gamma)+3$. 
\EndProof

\section{Proofs of the main results}\label{da proofs}

This section will include the proofs of Theorems A and B, which were
stated in the Introduction.

\Lemma\label{people who can't count}
Let $\Pi$ be a finitely generated, torsion-free subgroup of a group $\Pi'$.
Let $T$ be an element of $\Pi'$, 
and suppose that $\Pi'$ is
generated by $\Pi\cup\{T\}$ (so that $\Pi'$ is in particular
finitely generated). Let
$\barrho_0'$ be a representation of $\Pi'$ in $\pizzle(\CC)$ such that
$\barrho_0\doteq\barrho_0'|\Pi$ 
is a faithful representation of
$\Pi$. Suppose that $\barrho_0:\Pi\to\pizzle(\CC)$ admits a lift to
$\pizzle(\CC)$.
%Suppose that $\barrho_0$ and $\barrho_0'$ are smooth points of $\barR(\Pi)$
%and $\barR(\Pi')$ respectively, 
%\redcomment{
%I hope we don't need smoothness here. What we get from Lemma
%\ref{sorta marden} is that a certain kind of point of $\bar\barR(\Gamma)$
%lies in a unique irreducible component, and I hope that's what matters
%here. By the way, should this property of a point have a name?
%}
Let $V$ and $V'$ be
irreducible components of $\barR(\Pi)$
and $\barR(\Pi')$ containing
$\barrho_0$ and $\barrho_0'$ respectively, and assume that $V$ is the only
irreducible component of $\barR(\Pi)$ containing $\barrho_0$. Then either
% Suppose that somSuppose that some point of $V$ corresponds to a faithful
%representation $\barrho_0$ of $\Pi$, and that $\barrho_0$ admits an extension
%to a representation of $\Pi'$. Then either
\Alternatives
\item $\Pi'$ is the free product of the subgroups $\Pi$ and
  $\langle T\rangle$, or
\item  $\dim V'<(\dim V)+3$.
%the point of $R(\Pi')$ corresponding to  $\barrho_0'$  lies in an
%  irreducible 
%component of $R(\Pi)$ having dimension strictly less than $d+3$.
\EndAlternatives
\EndLemma

\Proof
Let $\langle T_0\rangle$ be an infinite cyclic group, and let $\Pi''$ denote the abstract free product $\Pi\star\langle
T_0\rangle$. Since $\Pi'$ is
generated by $\Pi\cup\{T\}$, there is a unique
surjective homomorphism from $h:\Pi''\to \Pi'$ which restricts to the inclusion homomorphism on
$\Pi$ and maps $T_0$ to $T$. 
Let $W$ denote the subset of $\barR(\Pi'')$ consisting of all
representations of $\Pi''$ whose restrictions to the kernel of $h$ are
trivial. Then $W$ is a Zariski-closed subset of $\barR(\Pi'')$, and there is an isomorphism of
algebraic sets $\alpha: \barR(\Pi')\to W$ defined by
$\alpha(\barrho)=\barrho\circ h$.
There is also an isomorphism of algebraic sets $\beta:\barR(\Pi'')\to
\barR(\Pi)\times\pizzle(\CC)$ defined by $\beta(\barrho)=(\barrho|\Pi,\barrho(T_0))$.
Thus $\beta\circ\alpha$ is an isomorphism of $\barR(\Pi')$ onto the
Zariski-closed subset $\beta( W)$ of $\barR(\Pi)\times\pizzle(\CC)$.
%Using the homomorphism $h$ we may identify
%$\Pi'$ with a quotient of $\Pi''$; then $\barR(\Pi')$
%becomes an algebraic subset of $\barR(\Pi'')$, consisting of all
%representations of $\Pi''$ whose kernels contain the kernel of $h$. 
%In turn, $\barR(\Pi'')$ may be identified with $\barR(\Pi)\times\pizzle(\CC)$: if $\barrho$
%is a representation of $\Pi$ in $\pizzle(\CC)$ and $A$ is an element of
%$\pizzle(\CC)$, the point $(\barrho,A)\in \barR(\Pi)\times\pizzle(\CC)$ is
%identified with the representation of $\Pi''=\Pi\star\langle
%T_0\rangle$ which restricts to $\barrho$ on $\Pi$ and maps $T_0$ to $A$. 
%Under these identifications we have
%$\barR(\Pi')\subset\barR(\Pi'')=\barR(\Pi)\times\pizzle(\CC)$; and
If we set
$E_0=\barrho_0'(T)$, we have 
%%\in
%%\barR(\Pi')\subset\barR(\Pi)\times\pizzle(\CC)$ 
%is equal to
%$
$\beta\circ\alpha(\barrho_0')=(\barrho_0'|\Pi, \barrho_0'(T))=
(\barrho_0,E_0)$.

Since $V$ is
the unique
irreducible component of $\barR(\Pi)$
containing
$\barrho_0$, the unique
irreducible component of $\barR(\Pi)\times\pizzle(\CC)$
containing
$(\barrho_0,E_0)$ is $V\times\pizzle(\CC)$. Since $V'$ is
an
irreducible component of $\barR(\Pi')
%\subset \barR(\Pi)\times\pizzle(\CC)
$
containing
$\barrho_0'$, 
the set
$\beta\circ\alpha(V')$ is
an
irreducible component of $\beta(W)\subset \barR(\Pi)\times\pizzle(\CC)$
containing
$\beta\circ\alpha(\barrho_0') =
(\barrho_0,E_0)$.
It follows that $\beta\circ\alpha(V')\subset V\times\pizzle(\CC)$.

Consider first the case in which $\beta\circ\alpha(V') $ is a proper subset of $
V\times\pizzle(\CC)$. In this case, since $V\times\pizzle(\CC)$ is
irreducible, we have $\dim V'=\dim \beta\circ\alpha(V')<\dim( V\times\pizzle(\CC))=(\dim V)+3$, so
that Alternative (ii) of the conclusion of the lemma
holds.

There remains the case in which $
\beta\circ\alpha(V')
=
V\times\pizzle(\CC)$. 
In particular we then have $\{\barrho_0\}
\times\pizzle(\CC)\subset 
\beta\circ\alpha(V')\subset \beta(W)$. 
Thus if for every
$E\in\pizzle(\CC)$ we set $\barsigma_E=\beta^{-1}(\barrho_0,E)$,
we have
 $\barsigma_E\in W$ for every
$E\in\pizzle(\CC)$. In view of the definition of $W$, it follows that:
\Claim\label{fuzzwonk}
$\ker\barsigma_E\supset\ker
h$ for every
$E\in\pizzle(\CC)$. 
\EndClaim
On the other hand, the definition of $\beta$ implies: 
\Claim\label{buzzsaw}
For every
$E\in\pizzle(\CC)$, the representation $\barsigma_E$ is the unique representation of
$\Pi''$ which restricts to $\barrho_0$  
on $\Pi$ and maps $T_0$ to $E$.
\EndClaim
%For every
%$E\in\pizzle(\CC)$, since $\barsigma_E\in W$, we have 
%$\barsigma_E\in \barR(\Pi')$, so
%that 

Now by hypothesis
%according to \cite[\redcomment{Get exact reference}]{splittings},
the  representation
$\barrho_0:\Pi\to\pizzle(\CC)$ admits a lift
$\rho_0:\Pi\to\zzle(\CC)$. Since $\barrho_0$ is faithful, $\rho_0$ is
also faithful, and $\rho_0(\Pi)$ contains no non-trivial scalar
matrix.
According to Lemma \ref{no but it will be}, there is an
element $L$ of
infinite order in $\zzle(\CC)$ such that the subgroup of $\zzle(\CC)$ generated by 
$\rho_0(\Pi)$ and $L$ is a free product of $\rho_0(\Pi)$ 
with the infinite cyclic
group $\langle L\rangle$, and contains no non-trivial scalar matrix.
If we denote by $E_1$ the image of $L$ under the quotient map from $\zzle(\CC)$ to
$\pizzle(\CC)$, it now follows that
the subgroup of $\pizzle(\CC)$ generated by  $\barrho_0(\Pi)$
and $E_1$ is a free product of $\barrho_0(\Pi)$ 
with the infinite cyclic
group $\langle E_1\rangle$. Applying \ref{buzzsaw} with $E=E_1$, we
deduce that $\barsigma_{E_1}$ is injective. 
But \ref{fuzzwonk}, applied with $E=E_1$, gives that 
$\ker
h\subset
\ker\barsigma_{E_1}$. Hence $h$ is injective in this case, so that Alternative (ii) of the conclusion of the lemma
holds.
\EndProof

\Proposition\label{agol bagol}
Let  $\Gamma$ be a finitely generated Kleinian group (so that
$\chibar(\Gamma)$ is defined by \ref{kleinian def}).
Then $\chibar(\Gamma)<\miof(\Gamma)$.
\EndProposition

\Proof
%\redcomment{
%The statement has changed. For the proof 
We first consider the case in
which $\Gamma$ has no rank-$2$ cusp subgroups.
%; in this case, the
%proof comes from fiddling with the argument below. 

 Set $k=\miof(\pi_1(N))$.
%Since $N$ is nearly simple, $\pi_1(N)$ is non-trivial, and is
%torsion-free by \ref{simple remarks}. Hence $k\ge1$ by
%\ref{well duh}. If $N$ is closed then $\chibar(N)=0$ and the
%conclusion follows. Hence we may assume that $N$ is not closed.
%Since $N$ is simple, it is a compact, aspherical triangulable space
%(see \ref{simple remarks}); hence $\pi_1(N)$ is homologically finite and
%$\chibar(N)=\chibar(\pi_1(N))$. 
By the definition of $\miof(\Gamma)$, there is a
finite generating set  $\Delta $  for $\Gamma$ such that $k=\iof(\Delta )$. In
particular, $\Delta $ contains $k$ independent elements. Hence we may write
$\Delta =\{x_1,\ldots,x_m\}$, where $m\ge k$ and $x_1,\ldots,x_k$ are
independent.

For $r=k,\ldots,m$, let $\Pi_r$ denote the subgroup $\langle
x_1,\ldots,x_r\rangle$ of $\Gamma$. We claim that
$\chibar(\Pi_r)<k$ for $r=k,\ldots,m$; the case $r=m$ of this claim is the
conclusion of the proposition, since $\Pi_m=\Gamma$. Note that $\Pi_k$ is free of
rank $k$, and hence $\chibar(\Pi_k)=k-1$; thus the claim is true
for $r=k$. It therefore suffices to prove that
$\chibar(\Pi_{r+1})\le\chibar(\Pi_r)$ whenever $k\le r<m$.

%Since $N$ is a nearly simple $3$-manifold and is not closed, it follows from
%Thurston's hyperbolization theorem \cite{morgan} that $N$ is homeomorphic to the compact
%core of $\HH^3/\Gamma_0$, where $\Gamma_0$ is a geometrically finite
%Kleinian group which has no parabolics and has infinite covolume. We shall fix
%an isomorphism $J: \pi_1(N)\to\Gamma_0$.

%It follows from \cite[Proposition 3.1.1]{splittings}
%that there is a subgroup $\hGamma_0$ of $\zzle(\CC)$ such that
%the quotient homomorphism $\zzle(\CC)\to\pizzle(\CC)$ restricts to an
%isomorphism $\beta:\hGamma_0\to\Gamma_0$.

Since the Kleinian group
$\Gamma$ has no rank-$2$ cusp subgroup, 
%is an infinite-covolume,
%geometrically finite Kleinian group, 
its finitely generated subgroups
$\Pi_r$  are also
Kleinian groups without rank-$2$ cusp subgroups. 
% geometrically finite by \cite[
%\redcomment{get precise reference}]{morgan}. They contain no parabolics since $\Gamma_0$ has
%none. 
It therefore follows from Lemma \ref{sorta marden} that, for any
integer $s$ with
$k\le s\le m$, if we
regard the inclusion homomorphism from 
$\Pi_s\le\Gamma$ to $\pizzle(\CC)$
%$\barrho_s=
%\beta^{-1}\circ 
%J|\Pi_s:\Pi_s\to\Gamma\le\pizzle(\CC)$
 as a faithful representation of $\Pi_s$ in
$\pizzle(\CC)$, then $\barrho_s$ lies in a unique irreducible
component $V_s$  of $\barR(\Pi_s)$, and that
\Equation\label{dim view}
\dim V_s=3\chibar(\Pi_s)+3.
\EndEquation

Suppose that $r$ is an integer with $k\le r<m$.  Since $\barrho_r$
is a discrete faithful representation, it follows from
\cite[Proposition 3.1.1]{splittings} that the representation $\barrho_r$ admits
a lift  to $\zzle(\CC)$. 
% we wish to apply Lemma \ref{people who can't 
The existence of such a lift, together with the fact that
$V_r$ is the only irreducible
component  of $\barR(\Pi_r)$ containing
$\barrho_r$, allows us to
apply Lemma \ref{people who can't
  count},
% which is required by the hypothesis
%of , 
%  count}, 
taking $\Pi=\Pi_r$, $\Pi'=\Pi_{r+1}$, $T=x_{r+1}$,
and letting $\barrho_r$ and $\barrho_{r+1}$ play the respective roles of
$\barrho_0$ and $\barrho'_0$. Lemma \ref{people who can't
  count} gives
%We may
%therefore apply Lemma \ref{people who can't
%  count}  
that either 
\Alternatives
\item $\Pi_{r+1}$ is the free product of the subgroups $\Pi_r$ and
  $\langle x_{r+1}\rangle$, or
\item  $\dim V_{r+1}<(\dim V_r)+3$.
\EndAlternatives
Thus for each $r$ with $k\le r<m$, either (i) or (ii) holds.

If (i) holds for a given $r$, then since $\Pi_k\le\Pi_r$
is generated by the independent elements $x_1,\ldots,x_k$, it follows
that $x_1,\ldots,x_{k+1}$ are independent. By the definition of index
of freedom, this means
that $\iof(\Delta )\ge k+1$, a contradiction to our choice of $\Delta
$. Hence for each $r$ with $k\le r<m$, Condition
(ii) must hold. Combining (ii) with the cases $s=r$ and $s=r+1$ of
(\ref{dim view}), we obtain
$3\chibar(\Pi_{r+1})+3<3\chibar(\Pi_r)+6$,
i.e. $\chibar(\Pi_{r+1})<\chibar(\Pi_r)+1$. Since $\chibar(\Pi_{r+1})$
and $\chibar(\Pi_r)$ are integers, this gives
$\chibar(\Pi_{r+1})\le\chibar(\Pi_r)$, as required. This completes the
proof of the proposition in the case in
which $\Gamma$ has no rank-$2$ cusp subgroups.

We now turn to the general case. Let $\Gamma$ be an arbitrary finitely
generated Kleinian group.
According to Lemma \ref{doodle poodle},
there exist
a finitely generated Kleinian
group $\Gamma_0$ with no rank-$2$ cusp subgroups, with
$\chibar(\Gamma_0)=\chibar(\Gamma)$, and a surjective homomorphism
$\eta:\Gamma\to\Gamma_0$.
Applying Lemma \ref{quotient miof}, with $G_1=\Gamma$ and
$G_1=\Gamma_0$, and defining $\eta$ as above, we deduce that
$\miof(\Gamma_0)\le\miof(\Gamma)$. Since 
$\Gamma_0$ has no rank-$2$ cusp subgroups, we may apply the
special case of the proposition that has already been proved to deduce
that $\chibar(\Gamma_0)<\miof(\Gamma_0)$. Thus we have
$\chibar(\Gamma)=\chibar(\Gamma_0) <\miof(\Gamma_0) \le\miof(\Gamma)$,
and the proof is complete.
\EndProof

We now prove Theorem B, which was stated in the introduction to this paper.

\Proof[Proof of Theorem B]
Since the finitely generated group $G$ is a subgroup of the
fundamental group of an orientable hyperbolic $3$-manifold, it is
isomorphic to a
Kleinian group. Thus the theorem
follows immediately from Proposition \ref{agol bagol}.
\EndProof

We are now in a position to prove Theorem A, which was stated in the introduction.

\Proof[Proof of Theorem A]
By Theorem B we have $\chibar(G)\le\miof(G)-1$. It therefore suffices
to show that $\miof(G)\le k-1$.

Assume that $\miof(G)\ge k$. Set $\Delta=\{[\alpha_1],\ldots,[\alpha_m]\}$. Since $\Delta$ is a generating set for $G$, we have $\iof(\Delta)\ge\miof(G)\ge k$. By definition this means that $\Delta$ contains $k$ independent elements. Thus we have $m\ge k$, and after possibly re-indexing the $\alpha_i$ we may assume that $[\alpha_1],\ldots,[\alpha_k]$ are independent.

Let us write $M=\HH^3/\Gamma$, where $\Gamma$ is a torsion-free
Kleinian group,  let $q:\HH^3\to M$ denote the quotient map, and let
us fix a point $z\in q^{-1}(p)$. Then $z$ defines an isomorphism
$j:\pi_1(M,p)\to\Gamma$. Set $\xi_i=j([\alpha_i])$ for
$i=1,\ldots,k$. Since $[\alpha_1],\ldots,[\alpha_k]$ are independent,
$\xi_1,\ldots,\xi_k$ freely generate a free subgroup of $\Gamma$, which may be regarded as a Kleinian group. If we set $d_i=\dist(z,\xi_i\cdot z)$ for $i=1,\ldots,k$, Theorem 4.1 of \cite{acs-surgery} now asserts that 
$$\sum_{i=1}^k\frac1{1+e^{d_i}}\le\frac12,$$
so that in particular $d_i\ge\log(2k-1)$ for some $i\in\{1,\ldots,k\}$. But for each $i\in\{1,\ldots,k\}$ we have $d_i\le\length\alpha_i$, which with the hypothesis of the theorem gives $d_i<\log(2k-1)$. This contradiction completes the proof.
\EndProof

\bibliographystyle{plain}

\begin{thebibliography}{10}

\bibitem{agol}
Ian Agol.
\newblock Tameness of hyperbolic 3-manifolds.
\newblock arXiv:math.GT/0405568.

\bibitem{acs-surgery}
Ian Agol, Marc Culler, and Peter~B. Shalen.
\newblock Dehn surgery, homology and hyperbolic volume.
\newblock {\em Algebr. Geom. Topol.}, 6:2297--2312, 2006.

\bibitem{accs}
James~W. Anderson, Richard~D. Canary, Marc Culler, and Peter~B. Shalen.
\newblock Free {K}leinian groups and volumes of hyperbolic {$3$}-manifolds.
\newblock {\em J. Differential Geom.}, 43(4):738--782, 1996.

\bibitem{threefree}
Gilbert Baumslag and Peter~B. Shalen.
\newblock Groups whose three-generator subgroups are free.
\newblock {\em Bull. Austral. Math. Soc.}, 40(2):163--174, 1989.

\bibitem{bznorms}
S.~Boyer and X.~Zhang.
\newblock On {C}uller-{S}halen seminorms and {D}ehn filling.
\newblock {\em Ann. of Math. (2)}, 148(3):737--801, 1998.

\bibitem{cg}
Danny Calegari and David Gabai.
\newblock Shrinkwrapping and the taming of hyperbolic 3-manifolds.
\newblock {\em J. Amer. Math. Soc.}, 19(2):385--446 (electronic), 2006.

\bibitem{splittings}
Marc Culler and Peter~B. Shalen.
\newblock Varieties of group representations and splittings of {$3$}-manifolds.
\newblock {\em Ann. of Math. (2)}, 117(1):109--146, 1983.

\bibitem{epstein}
D.~B.~A. Epstein.
\newblock Finite presentations of groups and {$3$}-manifolds.
\newblock {\em Quart. J. Math. Oxford Ser. (2)}, 12:205--212, 1961.




\bibitem{fps}
David Futer, Jessica~S. Purcell, and Saul Schleimer.
\newblock Effective drilling and filling of tame hyperbolic 3-manifolds.
\newblock {\em Comment. Math. Helv.}, 97(3):457--512, 2022.

\bibitem{ratioI}
Rosemary~K. Guzman and Peter~B. Shalen.
\newblock \title{The ratio of homology rank to hyperbolic volume, I}.
\newblock Journal of Topology and Analysis, To appear. arXiv:2207.00040.

\bibitem{ratioII}
Rosemary~K. Guzman and Peter~B. Shalen.
\newblock \title{The ratio of homology rank to hyperbolic volume, II: The role
  of the Four Color Theorem}.
\newblock Journal of Topology and Analysis, To appear. arXiv:2110.14847.

\bibitem{JS}
William~H. Jaco and Peter~B. Shalen.
\newblock Seifert fibered spaces in {$3$}-manifolds.
\newblock {\em Mem. Amer. Math. Soc.}, 21(220):viii+192, 1979.

\bibitem{jm}
Dennis Johnson and John~J. Millson.
\newblock Deformation spaces associated to compact hyperbolic manifolds.
\newblock In {\em Discrete groups in geometry and analysis ({N}ew {H}aven,
  {C}onn., 1984)}, volume~67 of {\em Progr. Math.}, pages 48--106.
  Birkh\"{a}user Boston, Boston, MA, 1987.

\bibitem{mishabook}
Michael Kapovich.
\newblock {\em Hyperbolic manifolds and discrete groups}.
\newblock Modern Birkh\"{a}user Classics. Birkh\"{a}user Boston, Ltd., Boston,
  MA, 2009.
\newblock Reprint of the 2001 edition.

\bibitem{namazi-souto}
Hossein Namazi and Juan Souto.
\newblock Non-realizability and ending laminations: proof of the density
  conjecture.
\newblock {\em Acta Math.}, 209(2):323--395, 2012.

\bibitem{newstead}
P.~E. Newstead.
\newblock {\em Introduction to moduli problems and orbit spaces}, volume~51 of
  {\em Tata Institute of Fundamental Research Lectures on Mathematics and
  Physics}.
\newblock Tata Institute of Fundamental Research, Bombay; Narosa Publishing
  House, New Delhi, 1978.

\bibitem{nisnevic}
V.~L. Nisnewitsch.
\newblock \"{U}ber {G}ruppen, die durch {M}atrizen \"uber einem kommutativen
  {F}eld isomorph darstellbar sind.
\newblock {\em Rec. Math. [Mat. Sbornik] N.S.}, 8 (50):395--403, 1940.

\bibitem{ohshika}
Ken'ichi Ohshika.
\newblock Realising end invariants by limits of minimally parabolic,
  geometrically finite groups.
\newblock {\em Geom. Topol.}, 15(2):827--890, 2011.

\bibitem{core}
G.~P. Scott.
\newblock Compact submanifolds of {$3$}-manifolds.
\newblock {\em J. London Math. Soc. (2)}, 7:246--250, 1973.

\bibitem{shalen-amalgam}
Peter~B. Shalen.
\newblock Linear representations of certain amalgamated products.
\newblock {\em J. Pure Appl. Algebra}, 15(2):187--197, 1979.

\bibitem{small}
Peter~B. Shalen.
\newblock Small optimal {M}argulis numbers force upper volume bounds.
\newblock {\em Trans. Amer. Math. Soc.}, 365(2):973--999, 2013.

\bibitem{wehrfritz}
B.~A.~F. Wehrfritz.
\newblock Generalized free products of linear groups.
\newblock {\em Proc. London Math. Soc. (3)}, 27:402--424, 1973.

\bibitem{wilson-pro-p}
John~S. Wilson.
\newblock On growth of groups with few relators.
\newblock {\em Bull. London Math. Soc.}, 36(1):1--2, 2004.

\bibitem{wilson-elem}
John~S. Wilson.
\newblock Free subgroups in groups with few relators.
\newblock {\em Enseign. Math. (2)}, 56(1-2):173--185, 2010.



\end{thebibliography}

\end{document}